\documentclass[a4paper,11pt]{article}
\usepackage{indentfirst,latexsym,bm,color}
\usepackage{amsmath,amssymb,amsfonts}

\usepackage[top=1in,left=1in,right=1in]{geometry}

\makeatletter

\@addtoreset{equation}{section} \makeatother

\begin{document}
\newtheorem{dingli}{Theorem}[section]
\newtheorem{dingyi}[dingli]{Definition}
\newtheorem{tuilun}[dingli]{Corollary}
\newtheorem{zhuyi}[dingli]{Remark}
\newtheorem{yinli}[dingli]{Lemma}

\title{  The structures and decompositions of symmetries involving idempotents
\thanks{ This work was supported by NSF of
China (Nos: 11671242, 11571211) and the
Fundamental Research Funds for the Central Universities
(GK201801011).}}
\author{Yuan Li$^a$\thanks{E-mail address:
 liyuan0401@aliyun.com}, \ Jiaxin Zhang$^a$, \ Nana Wei$^b$ }
\date{} \maketitle\begin{center}
\begin{minipage}{16cm}
{ \small $$ \ \ \ a. \ School \ of \ Mathematics \ and \
Information \ Science,\ Shaanxi \ Normal \ University, $$ $$ Xi'an,
710062, \ People's \ Republic \ of\ China.$$} {\small $$ \ \ \ b.  \ College \ of \ Xingzhi,\ Xi'an \ University \ of \ Finance \ and \ Economics,  $$ $$ Xi'an, \ 710038, \ People's \
Republic \ of \ China. $$ }
\end{minipage}
\end{center}
 \vspace{0.05cm}
\begin{center}
\begin{minipage}{16cm}
\ {\small {\bf Abstract }
Let $\mathcal{H}$ be a separable Hilbert space and $P$ be an idempotent on $\mathcal{H}.$ We denote by
$$\Gamma_{P}=\{J: J=J^{\ast}=J^{-1} \hbox{ }\hbox{ and }\hbox{ } JPJ=I-P\}$$
and $$\Delta_{P}=\{J: J=J^{\ast}=J^{-1} \hbox{ }\hbox{ and }\hbox{ } JPJ=I-P^*\}.$$
In this paper, we first get that symmetries $(2P-I)|2P-I|^{-1}$ and $(P+P^{*}-I)|P+P^{*}-I|^{-1}$ are the same. Then we show that $\Gamma_{P}\neq\emptyset$ if and only if $\Delta_{P}\neq\emptyset.$
Also, the specific structures of all symmetries
$J\in\Gamma_{P}$ and $J\in\Delta_{P} $ are established, respectively. Moreover,
 we prove that $J\in\Delta_{P}$ if and only if $\sqrt{-1}J(2P-I)|2P-I|^{-1}\in\Gamma_{P}.$
\\}
\endabstract
\end{minipage}\vspace{0.10cm}
\begin{minipage}{16cm}
\end{minipage}
\end{center}
\begin{center} \vspace{0.01cm}
\end{center}

\section{ Introduction}
Let $\mathcal{H}$ and $\mathcal{K}$ be separable complex Hilbert spaces,
and $\mathcal{B(H,K)}$
be the set of all bounded linear operators from
$\mathcal{H}$ into $\mathcal{K}.$ An operator $A\in\mathcal{B(H)}$ is called positive, if
$A\geqslant 0,$ meaning $\langle
Ax,x\rangle \geqslant0 $
 for all $x\in\mathcal{H},$ where $\langle , \rangle$ is the inner product of $\mathcal{H}.$
As usual, the operator order (Loewner partial order) relation $A\geqslant B$ between two self-adjoint operators is defined as $A-
B\geqslant0.$  Also, we denote by $\mathcal{B({H})}^{+}$ the set of all positive bounded linear operators on $\mathcal{H}.$  For $A \in\mathcal{B({H})}^{+},$
 $A^{\frac{1}{2}} $ is the positive square root of $A.$
 In particular, $|A|:=(A^{*}A)^{\frac{1}{2}}$ is the absolute value of operator $A,$
 where $A^{*}$ is the adjoint operator of $A.$
Let $S^+:=\frac{|S|+S}{2}$ and $S^-:=\frac{|S|-S}{2}$ be the positive and negative parts of a
self-adjoint operator $S\in\mathcal{B(H)}.$

For an operator $T\in \mathcal{B(H,K)}, N(T),R(T)$ and $\overline{R(T)}$ denote the null
 space, the range of $T,$ and the closure of $R(T),$ respectively.
we also use $P_{A}$ to denote the orthogonal projection onto
$\overline{R(A)}.$
 It is well known that every operator $T\in \mathcal{B(H,K)}$ has a (unique) polar decomposition
 $T = U(T^*T)^\frac{1}{2}$
where $U$ is a partial isometry from $H$ onto $\overline{R(T)}$ with kernel space $N(T).$
An operator $J\in \mathcal{B(H)}$ is said to be a symmetry (or self-adjoint unitary operator) if $J=J^*=J^{-1}.$
In this case, $J^+=\frac{I+J}{2}$ and $J^-=\frac{I-J}{2}$ are mutually annihilating orthogonal projections. If
$J$ is a non-scalar symmetry, then
 an indefinite inner product is defined by \[[x,y]:=\langle Jx, y\rangle \qquad(x,y\in \mathcal{H})\] and
 $(\mathcal{H}, J)$ is called a Krein space  \cite{s1,s2}.

Let $\mathcal{B(H)}^{Id}$ and $\mathcal{P(H)}$ be the set of all idempotents and orthogonal projections on $\mathcal{H},$ respectively. It is easy to see that the range of an
operator $P\in \mathcal{B(H)}^{Id}$ is closed and $P$ can be written as a $2\times2$ operator matrix:  \begin{equation}
P=\begin{bmatrix}\label{1.1}
I & P_1 \\
0 &  0\end{bmatrix}:R(P)\oplus R(P)^{\perp},\end{equation} where $P_1\in
{\mathcal{B}}(R(P)^{\perp}, R(P)).$ Also, we denote by
$$\Gamma_{P}:=\{J: J=J^{\ast}=J^{-1} \hbox{ }\hbox{ and }\hbox{ } JPJ=I-P\}$$
and $$\Delta_{P}:=\{J: J=J^{\ast}=J^{-1} \hbox{ }\hbox{ and }\hbox{ } JPJ=I-P^*\}.$$

In recent years, the descriptions for intertwining operators and a fixed difference properties
of two orthogonal projections are considered in \cite{s4,s7,s15,s16,s17}.
That is how to find a unitary operator $U$ such that $UP=QU$ and $UQ=PU$ for projections $P$ and $Q.$
For a pair $(P, Q)$ of orthogonal projections, the characterization of
intertwining operator $U$ is given in \cite{s4,s7,s14,s17}.
Indeed, these results also describe the sufficient and necessary condition for the existence of a symmetry $J$ with $JPJ=Q$ and the explicit formulas of all symmetries $J$ with $JPJ=Q.$
Furthermore, some decomposition properties of projections (or $J$-projections) were studied in \cite{s2,s6,s11}.
In particular, the existence of $J$-selfadjoint (positive, contractive) projections and its properties are obtained in \cite{s12,s13,s14}.
Also, the minimal and maximal elements of the set of all symmetries the symmetries $J$ with $P^{\ast}JP\leqslant J$ (or $JP\geqslant0)$ are given in \cite{s9,s10}.
 The purpose of this paper is to consider the structures and decompositions of a symmetry $J$ with $JPJ=I-P$ and $JPJ=I-P^*,$ respectively.
Firstly, we show that $(2P-I)|2P-I|^{-1}=(P+P^{*}-I)|P+P^{*}-I|^{-1},$
which means $\rho_{E}=s_{E}$ in \cite[Part 7]{s1}, so this equation is an extension of \cite[Proposition 7.1 and Remark 7.2]{s1}.
Then we show that $\Gamma_{P}\neq\emptyset$ if and only if $\Delta_{P}\neq\emptyset.$
Also, the specific structures of all symmetries
$J\in\Gamma_{P}$ and $J\in\Delta_{P} $ are characterized, respectively. Moreover,
 we get that $J\in\Delta_{P}$ if and only if $\sqrt{-1}J(2P-I)|2P-I|^{-1}\in\Gamma_{P}.$

\section{  The symmetries of $\Gamma_{P}$}

 In this section, we first get that $$(2P-I)|2P-I|^{-1}=(P+P^{*}-I)|P+P^{*}-I|^{-1}.$$
Then the structures of $J\in\Gamma_{P}$ are given.
At last, we present some equivalent conditions for the equation $UPU^*=Q,$ where $U$ is a unitary operator and $P,Q\in \mathcal{B(H)}^{Id}.$ To show our main results, the following two lemmas are needed.

{\bf Lemma 2.1.}  Let $P\in \mathcal{B(H)}^{Id}.$  Then $|2P^*-I|=|2P-I|^{-1}.$

{\bf Proof.}  Clearly, $(2P-I)^{2}=I.$  Then
 $(2P-I)^{-1}=2P-I,$ which implies
 $$|2P-I|^{-1}=[(2P^*-I)(2P-I)]^{-\frac{1}{2}}=[(2P-I)^{-1}(2P^*-I)^{-1}]^{\frac{1}{2}}=|2P^*-I|.$$
$\Box$

The following lemma is essentially from \cite[Proposition 3.1]{s5}. For the reader's convenience, we give a proof here.

{\bf Lemma 2.2.} Let $P\in \mathcal{B(H)}^{Id}.$  Then $J':=(2P-I)|2P-I|^{-1}$ is a symmetry and  $(2P-I)|2P-I|^{-1}=(2P^{*}-I)|2P^{*}-I|^{-1}.$

{\bf Proof.}
It is easy to see that  $$2P^*-I=|2P-I|^2(2P-I)=(2P-I)|2P^*-I|^2.$$ Thus $$J'^*=|2P-I|^{-1}(2P^*-I)=|2P-I|(2P-I)=J'^{-1}$$ and
$|2P-I|(2P-I)=(2P-I)|2P^*-I|,$ so $$J'^{-1}=|2P-I|(2P-I)=(2P-I)|2P^*-I|=(2P-I)|2P-I|^{-1}=J'.$$
 Then $J'$ is a symmetry. Also, $$(2P-I)|2P-I|^{-1}=(2P^{*}-I)|2P^{*}-I|^{-1}$$ follows from Lemma 2.1 and equation $2P-I=(2P^{*}-I)|2P-I|^{2}.$ $ \Box$

 The following result is an extension of [1, Remark 7.2].

{\bf Theorem 2.3.} Let $P\in \mathcal{B(H)}^{Id}.$
Then
 $(2P-I)|2P-I|^{-1}=(P+P^{*}-I)|P+P^{*}-I|^{-1}.$

{\bf Proof.}  It follows from [3] that $P+P^{*}-I$ is invertible.  Setting $J':=(2P-I)|2P-I|^{-1}$ and using Lemma 2.2, we get that
 $$2P-I=J'|2P-I|  \hbox{   }\hbox{   }  \hbox{ and }\hbox{   } \hbox{  }   2P^*-I=J'|2P^*-I|.$$ Then
\begin{equation}2(P+P^*-I)=J'(|2P-I|+|2P^*-I|).\end{equation}
Also Lemma 2.1 implies $$(\frac{|2P-I|+|2P^*-I|}{2})^2=\frac{|2P-I|^2+|2P^*-I|^2+2I}{4}=(P+P^{*}-I)^2,$$
so \begin{equation}\frac{|2P-I|+|2P^*-I|}{2}=|P+P^{*}-I|.\end{equation}
Combining equations (2.1) and (2.2), we have
$$P+P^*-I=J'|P+P^{*}-I|,$$ which yields  $$J'=(P+P^{*}-I)|P+P^{*}-I|^{-1}.$$

$\Box$

In the following, we give a concrete operator matrix form of the operator
$|2P-I|^{-1}.$

{\bf Corollary 2.4.} Let $P\in \mathcal{B(H)}^{Id}$ have the form (1.1).
Then \begin{equation}|2P-I|^{-1}=\left(\begin{array}{cc}
(I+2P_{1}P_{1}^{*})(I+P_{1}P_{1}^{*})^{-\frac{1}{2}}&-(I+P_{1}P_{1}^{*})^{-\frac{1}{2}}P_{1}\\-P_{1}^{*}(I+P_{1}P_{1}^{*})^{-\frac{1}{2}}&(I+P_{1}^{*}P_{1})^{-\frac{1}{2}}\end{array}\right):R(P)\oplus R(P)^{\perp}
\end{equation} and \begin{equation}|2P-I|=\left(\begin{array}{cc}
(I+P_{1}P_{1}^{*})^{-\frac{1}{2}}&(I+P_{1}P_{1}^{*})^{-\frac{1}{2}}P_{1}\\P_{1}^{*}(I+P_{1}P_{1}^{*})^{-\frac{1}{2}}&(I+2P_{1}^{*}P_{1})(I+P_{1}^{*}P_{1})^{-\frac{1}{2}}\end{array}\right):R(P)\oplus R(P)^{\perp}.
\end{equation}

{\bf Proof.} It is easy to verify that
$$
|P+P^{*}-I|=[\left(\begin{array}{cc}
I&P_{1}\\P_{1}^{*}&-I\end{array}\right)\left(\begin{array}{cc}
I&P_{1}\\P_{1}^{*}&-I\end{array}\right)]^{\frac{1}{2}}=
\left(\begin{array}{cc}
(I+P_{1}P_{1}^{*})^{\frac{1}{2}}&0\\0&(I+P_{1}^{*}P_{1})^{\frac{1}{2}}\end{array}\right),
$$
so
\begin{equation}
\begin{array}{rl}&(P+P^{*}-I)|P+P^{*}-I|^{-1}\\=&\left(\begin{array}{cc}
I&P_{1}\\P_{1}^{*}&-I\end{array}\right)\left(\begin{array}{cc}
(I+P_{1}P_{1}^{*})^{-\frac{1}{2}}&0\\0&(I+P_{1}^{*}P_{1})^{-\frac{1}{2}}\end{array}\right)\\=&
\left(\begin{array}{cc}
(I+P_{1}P_{1}^{*})^{-\frac{1}{2}}&P_{1}(I+P_{1}^{*}P_{1})^{-\frac{1}{2}}\\P_{1}^{*}(I+P_{1}P_{1}^{*})^{-\frac{1}{2}}&-(I+P_{1}^{*}P_{1})^{-\frac{1}{2}}\end{array}\right)
.\end{array}
\end{equation}
Clearly,
$
P_{1}(I+P_{1}^{*}P_{1})=(I+P_{1}P_{1}^{*})P_{1}.
$
Then
$
P_{1}(I+P_{1}^{*}P_{1})^{\frac{1}{2}}=(I+P_{1}P_{1}^{*})^{\frac{1}{2}}P_{1},
$
which induces,
\begin{equation}
P_{1}(I+P_{1}^{*}P_{1})^{-\frac{1}{2}}=(I+P_{1}P_{1}^{*})^{-\frac{1}{2}}P_{1}.
\end{equation}
 Thus Theorem 2.3 yields
 $$\begin{array}{rcl}|2P-I|^{-1}&=&(2P-I)(P+P^{*}-I)|P+P^{*}-I|^{-1}\\&=&\left(\begin{array}{cc}
(I+2P_{1}P_{1}^{*})(I+P_{1}P_{1}^{*})^{-\frac{1}{2}}&-(I+P_{1}P_{1}^{*})^{-\frac{1}{2}}P_{1}\\-P_{1}^{*}(I+P_{1}P_{1}^{*})^{-\frac{1}{2}}&(I+P_{1}^{*}P_{1})^{-\frac{1}{2}}\end{array}\right):R(P)\oplus R(P)^{\perp}.
\end{array}$$
 In a similar way, we get the formula of $|2P-I|.$ $\Box$

The following result is related to [9, Lemma 6] and [1, Proposition 5.1].

{\bf Proposition 2.5.} Let $P\in \mathcal{B(H)}^{Id}$ have the form (1.1). Then

(i) $N(P)\cap R(P)^\perp=N(P_{1})=N(P+P^{*}).$

(ii) $N(P)^\perp\cap R(P)=N(P_{1}^{*})=N(2I-P-P^{*}).$

{\bf Proof.} (i) Obviously, $N(P)\cap R(P)^\perp=N(P_1).$
It is easy to see that
$$ P+P^{*}=\left(\begin{array}{cc}
2I&P_{1}\\P_{1}^{*}&0\end{array}\right):R(P)\oplus R(P)^{\perp}.$$

Let $x\in R(P)$ and $y\in R(P)^{\bot}$ satisfy that
$$
\begin{array}{ll}(P+P^{*})\left(\begin{array}{c}
x\\y\end{array}\right)=\left(\begin{array}{cc}
2I&P_{1}\\P_{1}^{*}&0\end{array}\right)\left(\begin{array}{c}
x\\y\end{array}\right)=0
\end{array},
$$
which is $2x+P_{1}y=0$ and $P_{1}^{*}x=0$, so $2P_{1}^{*}x+P_{1}^{*}P_{1}y=0$, which implies $P_{1}^{*}P_{1}y=0.$
Thus
$ y \in N(P_{1}^{*}P_{1})=N(P_{1})$ i.e. $P_{1}y=0.$ Then $x=0,$ which says
 $N(P+P^{*})\subseteq0\oplus N(P_{1}).$
 Another inclusion relation $ 0\oplus N(P_{1})\subseteq N(P+P^{*})$ is clear.
 Thus $N(P+P^{*})=0\oplus N(P_{1})= N(P_{1}).$

(ii) Clearly,  equation (1.1) implies \begin{equation}\begin{array}{rl} I-P^*&=\left(\begin{array}{cc}0&0\\-P_1^*&I\end{array}\right):R(P)\oplus R(P)^{\perp}\\&=\left(\begin{array}{cc}I& -P_1^*\\0&0\end{array}\right):R(P)^{\perp}\oplus R(P)\\&=\left(\begin{array}{cc}I& -P_1^*\\0&0\end{array}\right):R(I-P^*)\oplus R(I-P^*)^{\perp},\end{array}\end{equation} since $R(P)^{\perp}=N(I-P)^{\perp}=R(I-P^*).$
Replacing $P$ of (i) by $I-P^*,$ we have
$$N(I-P^*)\cap R(I-P^*)^\perp=N(-P_{1}^*)=N(2I-P-P^{*}).$$ That is $N(P)^\perp\cap R(P)=N(P_{1}^{*})=N(2I-P-P^{*}).$
$\Box$

In what follows, some equivalent conditions for $\Gamma_{P}\neq\emptyset$ are given.

{\bf Theorem 2.6.} Let $P\in \mathcal{B(H)}^{Id}.$ Then the following statements are equivalent:

(i) There exits a symmetry $J$ such that $JPJ=I-P,$

(ii) $dimN(P+P^{*})=dimN(2I-P-P^{*}),$

(iii) There exits a unitary operator $U$ such that $UPU^{*}=I-P.$

{\bf Proof.} $(i)\Rightarrow(iii)$ is obvious.

$(iii)\Rightarrow(ii).$ If $UPU^{*}=I-P,$ then
$UP^{*}U^{*}=(UPU^{*})^{*}=(I-P)^{*}=I-P^{*},$
so
$$U(P+P^{*})U^{*}=2I-P-P^{*}.$$

If $x\in N(2I-P-P^{*})$, then $U(P+P^{*})U^{*}x=0$, which yields $(P+P^{*})U^{*}x=0,$ so $U^{*}x\in N(P+P^{*}).$ Thus $$ dimN(2I-P-P^{*})\leq dimN(P+P^{*})$$ follows from that $U$ is a unitary operator.

Similarly, $P+P^{*}=U^{*}(2I-P-P^{*})U$ implies $ dimN(P+P^{*})\leq dimN(2I-P-P^{*}),$ so (ii) holds.

$(ii)\Rightarrow(i).$ Suppose that $P\in \mathcal{B(H)}^{Id}$ has the form (1.1).
By Proposition 2.5, $dimN(P+P^{*})=dimN(2I-P-P^{*})$ implies $dimN(P_{1})=dimN(P_{1}^{*}),$
so there exists a unitary operator $V$ such that $P_{1}^{*}=V(P_{1}P_{1}^{*})^{\frac{1}{2}}.$

Let $U:=V^{*}$. Then $UP_{1}^{*}=(P_{1}P_{1}^{*})^{\frac{1}{2}}=P_{1}U^{*},$ which yields $$UP_{1}^{*}P_{1}=P_{1}U^{*}P_{1}=P_{1}P_{1}^{*}U,$$ so
\begin{equation}U(I+P_{1}^{*}P_{1})^{-\frac{1}{2}}=(I+P_{1}P_{1}^{*})^{-\frac{1}{2}}U.\end{equation}
We define an operator $J_{11}$ from $R(P)$ into $R(P)$ by $$J_{11}:=-U(I+P_{1}^{*}P_{1})^{-\frac{1}{2}}P_{1}^{*}.$$
Then
\begin{align*}J_{11}&=-U(I+P_{1}^{*}P_{1})^{-\frac{1}{2}}P_{1}^{*}\\&=
-(I+P_{1}P_{1}^{*})^{-\frac{1}{2}}UP_{1}^{*}\\&=
-(I+P_{1}P_{1}^{*})^{-\frac{1}{2}}P_{1}U^{*}\\&=
-P_{1}(I+P_{1}^{*}P_{1})^{-\frac{1}{2}}U^{*}\\&=J_{11}^{*} ,\end{align*}
so $J_{11}$ is self-adjoint.

Let $J_{12}$ from $R(P)$ into $R(P)^\perp$ have the form
 $$J_{12}:=U(I+P_{1}^{*}P_{1})^{-\frac{1}{2}}$$
and $J_{22}$ from $R(P)^\perp$ into $R(P)^\perp$ have the form
$$J_{22}:=P_{1}^{*}U(I+P_{1}^{*}P_{1})^{-\frac{1}{2}},$$ respectively.
In a similar way, we get that
\begin{align*}
J_{22}&=P_{1}^{*}U(I+P_{1}^{*}P_{1})^{-\frac{1}{2}}\\&=
U^{*}P_{1}(I+P_{1}^{*}P_{1})^{-\frac{1}{2}}\\&=
U^{*}(I+P_{1}P_{1}^{*})^{-\frac{1}{2}}P_{1}\\&=
(I+P_{1}^{*}P_{1})^{-\frac{1}{2}}U^{*}P_{1}\\&=
J_{22}^{*} ,
\end{align*}
so $J_{22}$ is self-adjoint.
Define a self-adjoint operator $J\in \mathcal{B(H)}$ by
$$ J=\left(\begin{array}{cc}
J_{11}&J_{12}\\J_{12}^{*}&J_{22}\end{array}\right)=\left(\begin{array}{cc}
-U(I+P_{1}^{*}P_{1})^{-\frac{1}{2}}P_{1}^{*}&U(I+P_{1}^{*}P_{1})^{-\frac{1}{2}}\\(I+P_{1}^{*}P_{1})^{-\frac{1}{2}}U^{*}&P_{1}^{*}U(I+P_{1}^{*}P_{1})^{-\frac{1}{2}}\end{array}\right)
:R(P)\oplus R(P)^{\perp}.$$
Then a direct calculation yields $J^{2}=I$ and $JPJ=(I-P).$
\ \ \ $\Box$

In the following, we shall get the specific operator matrix forms for symmetries $J\in\Gamma_{P}.$

{\bf Theorem 2.7.} Let $P\in \mathcal{B(H)}^{Id}$ have the form (1.1).
If $JPJ=I-P$ for some symmetries $J,$ then $dimR(P)=dimR(P)^{\bot}$ and
$$J=\left(\begin{array}{cc}
-U(I+P_{1}^{*}P_{1})^{-\frac{1}{2}}P_{1}^{*}&U(I+P_{1}^{*}P_{1})^{-\frac{1}{2}}
\\(I+P_{1}^{*}P_{1})^{-\frac{1}{2}}U^{*}&P_{1}^{*}U(I+P_{1}^{*}P_{1})^{-\frac{1}{2}}\end{array}\right):R(P)\oplus R(P)^{\perp}
,$$
where $U$ is a unitary operator from $R(P)^{\bot}$ onto $R(P)$ with $U^{*}P_{1}=P_{1}^{*}U.$

{\bf Proof.}  Suppose that $J$ has the following operator matrix form  $$J=\left(\begin{array}{cc}
J_{11}&J_{12}\\J_{12}^{*}&J_{22}\end{array}\right):R(P)\oplus R(P)^{\perp},$$
where $J_{11}$ and $J_{22}$ are self-adjoint operators.
 It follows from the fact $JP=(I-P)J$ that
\begin{equation}\begin{cases}J_{11}=-P_{1}J_{12}^{*}\qquad\qquad\ \ \ \  \textcircled{1}
&\\J_{11}P_{1}=-P_{1}J_{22}\;\qquad \qquad\textcircled{2}&\\J_{12}^{\ast}P_{1}=J_{22}\ \ \ \ \ \ \qquad\qquad \textcircled{3}.\end{cases}\end{equation}

On the other hand, $J=J^{\ast}=J^{-1}$ yields $J^{2}=I,$  so
 \begin{equation}\begin{cases}J_{11}^2+J_{12}J_{12}^{\ast}=I\qquad\qquad \textcircled{1}
&\\J_{11}J_{12}+J_{12}J_{22}=0\;\ \qquad \textcircled{2}&\\J_{12}^{\ast}J_{12}+J_{22}^2=I\qquad\qquad \textcircled{3}.\end{cases}\end{equation}

It is easy to see that equations $\textcircled{1}$ of (2.9) and $\textcircled{1}$ of (2.10) imply
$$
J_{12}(P_{1}^{*}P_{1}+I)J_{12}^{*}=J_{12}P_{1}^{*}P_{1}J_{12}^{*}+J_{12}J_{12}^{*}=J_{11}^2+J_{12}J_{12}^{\ast}=I,
$$
so $J_{12}$ is right invertible.

Using equations $\textcircled{1},$ $\textcircled{2}$ and $\textcircled{3}$ of (2.9), we have
that \begin{equation}J_{12}P_{1}^{*}P_{1}=-J_{11}P_1=P_1J_{22}=P_1J_{22}^*=P_{1}P_{1}^{*}J_{12},\end{equation} which induces
$$
J_{12}(P_{1}^{*}P_{1}+I)J_{12}^{*}=(I+ P_{1}P_{1}^{*})J_{12}J_{12}^{*}=I.
$$
Thus
\begin{equation}
J_{12}J_{12}^{*}=(I+ P_{1}P_{1}^{*})^{-1}.
\end{equation}

Moreover, by equations $\textcircled{3}$ of (2.9) and $\textcircled{3}$ of (2.10), we know that  \begin{equation}
J_{12}^{*}(I+P_{1}P_{1}^{*})J_{12}=J_{12}^{*}J_{12}+J_{12}^{*}P_{1}P_{1}^{*}J_{12}
=J_{12}^{*}J_{12}+J_{22}^2=I,
\end{equation}
so $J_{12}$ is left invertible. Then $J_{12}$ is invertible, which yields $dimR(P)=dimR(P)^{\bot}.$  Also, equations (2.11) and (2.13) imply
$$
(I+P_{1}^{*}P_{1})J_{12}^{*}J_{12}=J_{12}^{*}(I+P_{1}P_{1}^{*})J_{12}=I,
$$
which means
\begin{equation}
J_{12}^{*}J_{12}=(I+P_{1}^{*}P_{1})^{-1}.
\end{equation}
Using equation (2.14) and the polar decomposition theorem,
we conclude that there exits a uniquely unitary operator $U$ such that
\begin{equation}
J_{12}=U(J_{12}^{*}J_{12})^{\frac{1}{2}}=U(I+P_{1}^{*}P_{1})^{-\frac{1}{2}}.
\end{equation}
Then from equations $\textcircled{1}$ and $\textcircled{3}$ of (2.9), we get that
$$
J_{11}=-J_{12}P_{1}^{*}=-U(I+P_{1}^{*}P_{1})^{-\frac{1}{2}}P_{1}^{*}
 \hbox{   }\hbox{  }\hbox{  }\hbox{ and }\hbox{  }\hbox{  }\hbox{  }
J_{22}=
P_{1}^{*}J_{12}=P_{1}^{*}U(I+P_{1}^{*}P_{1})^{-\frac{1}{2}},
$$
which implies
\begin{equation}
J_{22}=J_{22}^*=(I+P_{1}^{*}P_{1})^{-\frac{1}{2}}U^{*}P_{1}.
\end{equation} Furthermore, \begin{equation}J_{22}=P_{1}^{*}U(I+P_{1}^{*}P_{1})^{-\frac{1}{2}}=P_{1}^{*}(I+P_{1}P_{1}^{*})^{-\frac{1}{2}}U
=(I+P_{1}^{*}P_{1})^{-\frac{1}{2}}P_{1}^{*}U,
\end{equation} follows from equations (2.12) and (2.15).
Thus $U^{*}P_{1}=P_{1}^{*}U.$
\ \ \ $\Box$

{\bf Corollary 2.8.} Let $P\in \mathcal{B(H)}^{Id}$ have the form (1.1). If $ dimN(P+P^{*})=dimN(2I-P-P^{*})=0$, then $JPJ=I-P$ if and only if
$$ J=\left(\begin{array}{cc}
-VJ_{1}(I+P_{1}^{*}P_{1})^{-\frac{1}{2}}P_{1}^{*}&VJ_{1}(I+P_{1}^{*}P_{1})^{-\frac{1}{2}}
\\(I+P_{1}^{*}P_{1})^{-\frac{1}{2}}J_{1}V^{*}&P_{1}^{*}VJ_{1}(I+P_{1}^{*}P_{1})^{-\frac{1}{2}}\end{array}\right):R(P)\oplus R(P)^{\perp}
,$$
where $V$ is the unique unitary operator with $P_{1}=V(P_{1}^{*}P_{1})^{-\frac{1}{2}}$ and
$J_{1}\in\mathcal{B}( R(P)^{\perp})$ is a symmetry with $J_{1}P_{1}^{*}P_{1}=P_{1}^{*}P_{1}J_{1}.$

{\bf Proof.} Necessity.  If $U\in \mathcal{B}(R(P)^{\bot},R(P))$ is a unitary operator  with $P_{1}^{*}U=U^{*}P_{1},$ then \begin{equation}P_{1}^{*}P_{1}=P_{1}^{*}UU^{*}P_{1}=U^{*}P_{1}P_{1}^{*}U.\end{equation}
Using Proposition 2.5, we conclude that $P_{1}$ is a injective and dense range.
Thus the polar decomposition theorem implies $P_{1}=V(P_{1}^{*}P_{1})^{\frac{1}{2}},$ where $V$ is the unique  unitary, which yields
 \begin{equation}P_{1}P_{1}^{*}=V(P_{1}^{*}P_{1})^{\frac{1}{2}}(P_{1}^{*}P_{1})^{\frac{1}{2}}V^{*}
 =V(P_{1}^{*}P_{1})V^{*}.\end{equation}
Then equations (2.18) and (2.19) imply $$
P_{1}^{*}P_{1}=U^{*}P_{1}P_{1}^{*}U=U^{*}V(P_{1}^{*}P_{1})V^{*}U,
$$
so
$$
(P_{1}^{*}P_{1})( V^{*}U)=( V^{*}U) (P_{1}^{*}P_{1}).
$$

Furthermore, equation $P_{1}^{*}U=U^{*}P_{1}$ also yields
$(P_{1}^{*}P_{1})^{\frac{1}{2}}V^{*}U=U^{*}V(P_{1}^{*}P_{1})^{\frac{1}{2}},$ which induces
$$
U^{*}V(P_{1}^{*}P_{1})=(P_{1}^{*}P_{1})V^{*}U=V^{*}U(P_{1}^{*}P_{1}),
$$
so  $U^{*}V=V^{*}U$ follows from the fact that $P_1^*$ is a dense range.
Setting $J_{1}:=U^{*}V,$ we know that $J_{1}$ is a symmetry with
$U=VJ_{1} $ and $(P_{1}^{*}P_{1})J_1=J_1(P_{1}^{*}P_{1}).$
Then the necessity follows from Theorem 2.7.

Sufficiency. If $J'\in\mathcal{B}( R(P)^{\perp})$ is a symmetry with $J'P_{1}^{*}P_{1}=P_{1}^{*}P_{1}J'$ and $U=VJ',$ then
$$P_{1}^{*}U=P_{1}^{*}VJ'=(P_{1}^{*}P_{1})^{\frac{1}{2}}J'=J'(P_{1}^{*}P_{1})^{\frac{1}{2}}=J'V^{*}V(P_{1}^{*}P_{1})^{\frac{1}{2}}
=U^{*}P_{1},$$ so the proof of $(ii)\Longrightarrow (i)$ in Theorem 2.6 implies the sufficiency.
\ \ \ $\Box$

 At the last of this section, we give an equivalent condition for the equation $UPU^*=Q,$ where $U$ is a unitary operator and $P,Q\in \mathcal{B(H)}^{Id}.$ We need the following well-known lemma.

{\bf Lemma 2.9.}([1,3]) Let $P\in \mathcal{B(H)}^{Id}.$ Then

(i) $P_{N(P)}-P_{R(P)}$ and $P+P^*-I$ are invertible with $(P+P^*-I)^{-1}=P_{R(P)}-P_{N(P)}.$

(ii) $P=P_{R(P)}(P_{R(P)}-P_{N(P)})^{-1}$ and $P_{N(P)}=-(I-P)(P+P^*-I)^{-1}.$

{\bf Proposition 2.10.} Let $P,Q\in \mathcal{B(H)}^{Id}$ and $U$ be a unitary operator. Then the following statements are equivalent:

(i) $UPU^*=Q,$

(ii) $UP_{R(P)}U^*=P_{R(Q)}$ and $UP_{N(P)}U^*=P_{N(Q)},$

(iii) $U(P_{R(P)}-P_{N(P)})U^*=P_{R(Q)}-P_{N(Q)}$ and $U(P-P^{*})U^*=Q-Q^{*}.$

{\bf Proof.} $(ii)\Rightarrow(i).$ If $UP_{R(P)}U^*=P_{R(Q)}$ and $UP_{N(P)}U^*=P_{N(Q)},$ then $$U(P_{R(P)}-P_{N(P)})U^*=(P_{R(Q)}-P_{N(Q)}),$$ so
$$U(P_{R(P)}-P_{N(P)})^{-1}U^*=(P_{R(Q)}-P_{N(Q)})^{-1},$$ follows from Lemma 2.9 (i).
 Also, Lemma 2.9 (ii) implies \begin{align*}
UPU^*&=UP_{R(P)}(P_{R(P)}-P_{N(P)})^{-1}U^*
=(UP_{R(P)}U^*)(U(P_{R(P)}-P_{N(P)})^{-1}U^*)
\\&=P_{R(Q)}(P_{R(Q)}-P_{N(Q)})^{-1}=
Q.
\end{align*}

$(i)\Rightarrow(ii).$ If $UPU^*=Q,$ then $UP^{*}U^*=Q^{*},$ so $$U(P+P^{*}-I)U^*
=Q+Q^{*}-I,$$ which implies
 $$U(P+P^{*}-I)^{-1}U^*=(Q+Q^{*}-I)^{-1}.$$
Using Lemma 2.9 again, we get that
\begin{align*}
UP_{R(P)}U^*&=UP(P+P^{*}-I)^{-1}
U^*=(UPU^*)[U(P+P^{*}-I)^{-1}U^*]
\\&=Q(Q^{*}+Q-I)^{-1}\\&= P_{R(Q)}.
\end{align*} Similarly, we also have $UP_{N(P)}U^*=P_{N(Q)}.$
If $(i)$ and $(ii)$ hold, then $(iii)$ is obvious.

$(iii)\Rightarrow(i).$ If $U(P_{R(P)}-P_{N(P)})U^*=P_{R(Q)}-P_{N(Q)},$ then $$U(P+P^{*}-I)U^*=U(P_{R(P)}-P_{N(P)})^{-1}U^*
=(P_{R(Q)}-P_{N(Q)})^{-1}=Q+Q^{*}-I$$
follows from Lemma 2.9 (i). Thus $$U(P+P^{*})U^*=Q+Q^{*}.$$
Combining with the fact of $ U(P-P^{*})U^{*}=Q-Q^{*},$  we get that
$UPU^*=Q$ as desired.
\ \ \ $\Box$

\section{ The symmetries of $\Delta_{P}$ }

In this section, the structures of symmetries  $J\in\Delta_{P}$ are considered.
Moreover, we present some decomposition properties between symmetries of $\Gamma_{P}$
 and symmetries of $\Delta_{P}.$
  The following two lemmas are needed.

{\bf Lemma 3.1.} Let $A\in \mathcal{B(H,K)}.$
Then there exists a unitary operator $U\in \mathcal{B(K,H)}$
such that $UA$ is self-adjoint if and only if $dim N(A)=dim N(A^*).$

{\bf Proof.}   Sufficiency is clear from the polar decomposition theorem [8, Problem 134].

Necessity. If $UA$ is self-adjoint, then $UA=A^*U^*,$ so $UAU=A^*.$
Let $x\in N(A^*).$  Obviously, $UAUx=0,$
which yields $AUx=0,$ so $Ux\in N(A).$
Thus  $dim N(A)\geq dim N(A^*).$
Similarly, we get that $dim N(A)\leq dim N(A^*)$ from equation $A=U^*A^*U^*.$ Thus $dim N(A)=dim N(A^*).$ $\Box$

{\bf Lemma 3.2.} Let $A\in \mathcal{B(H,K)}$ and $U\in \mathcal{B(K,H)}$ be
 a unitary operator.
 Then $UA$ is self-adjoint if and only if $U=J_1V_1^*\oplus U_2$
 with respect to the space decomposition $\mathcal{K}=\overline{R(A)}\oplus R(A)^{\perp}$ and
 $\mathcal{H}=N(A)^{\perp}\oplus N(A),$ where $V_1$ is
 the unique unitary operator from
 $\overline{R(A^*)}$ into $\overline{R(A)}$ with $A=V_1(A^*A)^\frac{1}{2},$
 $J_1\in \mathcal{B}(N(A)^{\perp})$ is a symmetry with $(J_1\oplus 0)A^*A=A^*A(J_1\oplus 0)$ and  $U_2$ is
a unitary operator from
 $R(A)^{\perp}$ into $N(A).$

{\bf Proof.} Sufficiency is clear.

Necessity. It is clear that $A$ has the operator matrix form
$$A=\left(\begin{array}{cc}
A_1&0\\0&0\end{array}\right):N(A)^{\perp}\oplus N(A)\rightarrow\overline{R(A)}\oplus R(A)^{\perp}$$
and $ A_1\in B(N(B)^{\perp},\overline{R(B)})$ is injective and dense range.
 Suppose that $U\in \mathcal{B(K,H)}$ has the operator matrix form
$$U=\left(\begin{array}{cc}
U_{11}&U_{12}\\U_{21}&U_{22}\end{array}\right):\overline{R(A)}\oplus R(A)^{\perp}\rightarrow N(A)^{\perp}\oplus N(A).$$
By a direct calculation, we know that equation
$UA=A^*U^*$ implies that $U_{11}A_1=A_1^*U_{11}^*$ and $U_{21}A_1=0,$ so $U_{21}=0.$

 By Lemma 3.1, we have $dim N(A)=dim N(A^*).$  Thus the polar decomposition theorem implies that there exists  a unitary operator $V=V_1\oplus V_2$ such that $A=V(A^*A)^\frac{1}{2},$ where $V_1$ is
 the unique unitary operator from
 $\overline{R(A^*)}$ onto $\overline{R(A)}$ and $V_2$ is
 any unitary operator from
 $N(A)$ onto $N(A^*).$
 Then $AA^*=V(A^*A)V^*.$ Also, $UAA^*U^*=A^*A$ follows from $UA=A^*U^*.$
 Thus $UV(A^*A)V^*U^*=A^*A,$ which yields $UV(A^*A)=A^*AUV.$
On the other hand, $$UV(A^*A)^\frac{1}{2}=UA=A^*U^*=(A^*A)^\frac{1}{2}V^*U^*,$$
which implies $UV(A^*A)=(A^*A)V^*U^*,$ so $A^*A(UV-V^*U^*)=0.$
Obviously, $$UV-V^*U^*=\left(\begin{array}{cc}
U_{11}V_1-V_1^*U_{11}^*&U_{12}V_2\\-V_2^* U_{12}^*&U_{22}V_2-V_2^*U_{22}^*\end{array}\right):N(A)^{\perp}\oplus N(A)\rightarrow N(A)^{\perp}\oplus N(A),$$
which yields that $$A_1(U_{11}V_1-V_1^*U_{11}^*)=0
\hbox{   }\hbox{ and }\hbox{   } A_1U_{12}V_2=0,$$ from equation $A(UV-V^*U^*)=0.$
Thus $U_{11}V_1-V_1^*U_{11}^*=0$ and $U_{12}=0,$ so $U=U_{11}\oplus U_{22},$
which induces that $U_{11}$ and $U_{22}$ are unitary operators.
Setting $U_2=U_{22}$ and $J_1:=U_{11}V_1,$ we get that $J_1=J_1^*=J_1^{-1}$ and $U_{11}=J_1V_1^*.$
It is easy to see that $$J_1(A_1^*A_1)^\frac{1}{2}=U_{11}V_1(A_1^*A_1)^\frac{1}{2}
=U_{11}A_1=A_1^*U_{11}^*=A_1^*V_1V_1^*U_{11}^*=(A_1^*A_1)^\frac{1}{2}J_1,$$
which implies $(J_1\oplus 0)A^*A=A^*A(J_1\oplus 0).$   $\Box$

{\bf Proposition 3.3.}   Let $A\in \mathcal{B(H,K)}$ and $U\in \mathcal{B(K,H)}$ be
 a unitary operator.
 Then $UA\geq 0$ if and only if $U=V_1^*\oplus U_2$
 with respect to the space decomposition $\mathcal{K}=\overline{R(A)}\oplus R(A)^{\perp}$ and
 $\mathcal{H}=N(A)^{\perp}\oplus N(A),$ where $V_1$ is
 the unique unitary operator from
 $\overline{R(A^*)}$ into $\overline{R(A)}$ with $A=V_1(A^*A)^\frac{1}{2}$  and  $U_2$ is
a unitary operator from
 $R(A)^{\perp}$ into $N(A).$

{\bf Proof.} Sufficiency is clear.

Necessity.  By Lemma 3.2, we only need to show $J_1=I_1,$  where $I_1$
is the identity operator of subspace $N(A)^{\perp}.$
Using Lemma 3.2 again, we have $$J_1(A_1^*A_1)^\frac{1}{2}=(A_1^*A_1)^\frac{1}{2}J_1=A_1^*U_{11}^*\geq 0,$$      so $J_1\geq0$ follows from the fact that $A_1$ is injective.
 Thus $J_1=I_1,$ as $J_1\in \mathcal{B}(N(A)^{\perp})$ is a symmetry. $\Box$

 Some equivalent conditions for $\Delta_{P}\neq\emptyset$ are given in the following theorem.
It is worth noting that
$\Delta_{P}\neq\emptyset$  if and only if $\Gamma_{P}\neq\emptyset.$

{\bf Theorem 3.4.} Let $P\in \mathcal{B(H)}^{Id}.$ Then the following statements are equivalent:

(i) There exits a symmetry $J$ such that $JPJ=I-P^*,$

(ii) $dim N(P+P^{*})=dim N(2I-P-P^{*}),$

(iii) There exits a unitary operator $U$ such that $UPU^{*}=I-P^*,$

{\bf Proof.} $(i)\Rightarrow(iii)$ is obvious.

$(iii)\Rightarrow(ii)$ is analogous to that of Theorem 2.6.

$(ii)\Rightarrow(i).$ Suppose that $P\in \mathcal{B(H)}^{Id}$ has the form (1.1). By Proposition 2.5,
$dimN(P+P^{*})=dimN(2I-P-P^{*})$ implies $dimN(P_{1})=dimN(P_{1}^{*}),$ so there exists a unitary operator $V$
such that $P_{1}^{*}=V(P_{1}P_{1}^{*})^{\frac{1}{2}}$ and $V=V_1\oplus V_2,$ where $V_1$ is
 the unique unitary operator from
 $\overline{R(P_{1})}$ into $\overline{R(P_{1}^*)}$  and  $V_2$ is
a unitary operator from
 $R(P_{1})^{\perp}$ into $N(P_{1}).$
 Defining a self-adjoint operator $J$ by $$J:=\left(\begin{array}{cc}
0&\sqrt{-1}V^*
\\ -\sqrt{-1}V&0\end{array}\right):R(P)\oplus
R(P)^{\perp},$$ we get that $J^2=I$ and $$JPJ=\left(\begin{array}{cc}
0&0
\\ -VP_1V& I\end{array}\right)=\left(\begin{array}{cc}
0&0
\\ -P_1^*& I\end{array}\right)=I-P^*,$$ since
$V^*P_1^*=(P_{1}P_{1}^{*})^{\frac{1}{2}}=P_1V$ implies $ VP_1V=P_1^*.$
$\Box$

In the
following, we shall get a concrete operator matrix form for symmetries $J\in\Delta_{P}.$

{\bf Theorem 3.5.} Let $P\in \mathcal{B(H)}^{Id}$ have the form (1.1) and $J$ be a symmetry.
If $JPJ=I-P^*,$ then
$dimR(P)=dimR(P)^{\bot}$ and $$J=\left(\begin{array}{cc}
0&\sqrt{-1}U^*
\\ -\sqrt{-1}U&0\end{array}\right):R(P)\oplus
R(P)^{\perp},$$  where $U=J_1V_1^*\oplus U_2$
 with respect to the space decomposition $R(P)=\overline{R(P_1)}\oplus R(P_1)^{\perp}$ and
 $R(P)^{\perp}=N(P_1)^{\perp}\oplus N(P_1),$ where $V_1$ is
 the unique unitary operator from
 $\overline{R(P_1^*)}$ into $\overline{R(P_1)}$ with $P_1=V_1(P_1^*P_1)^\frac{1}{2},$
 $J_1\in \mathcal{B}(N(P_1)^{\perp})$ is a symmetry with $(J_1\oplus 0)P_1^*P_1=P_1^*P_1(J_1\oplus 0)$ and  $U_2$ is
a unitary operator from
 $R(P_1)^{\perp}$ into $N(P_1).$

{\bf Proof.}  Suppose that $J$ has the following operator matrix form  $$J=\left(\begin{array}{cc}
J_{11}&J_{12}\\J_{12}^{*}&J_{22}\end{array}\right):R(P)\oplus R(P)^{\perp},$$ where $J_{11}$ and $J_{22}$ are
self-adjoint operators.
 It follows from equation $JP=(I-P^*)J$ that \begin{equation} \left(\begin{array}{cc}
J_{11}&J_{11}P_1
\\ J_{12}^*&J_{12}^*P_1\end{array}\right)=\left(\begin{array}{cc}
0&0
\\ -P_1^*J_{11}+J_{12}^*&-P_1^*J_{12}+J_{22}\end{array}\right),\end{equation}
 which yields $J_{11}=0$ and $J_{12}^*P_1+P_1^*J_{12}=J_{22}.$

On the other hand,  $J^{2}=I$  implies
 \begin{equation}\begin{cases}J_{12}J_{12}^{\ast}=I \ \ \qquad\qquad\qquad \textcircled{1}
&\\J_{12}J_{22}=0\;\ \qquad \qquad\qquad\textcircled{2}&\\J_{12}^{\ast}J_{12}+J_{22}^2=I\qquad\qquad
\textcircled{3}.\end{cases}\end{equation}

From equations $\textcircled{1}$ of (3.2),  we have $J_{12}$ is co-isometry,
so $J_{12}^*J_{12}$ is an orthogonal projection.
 Thus $J_{22}^2=I-J_{12}^*J_{12}$ is also an orthogonal projection.
 As $J_{22}$ is self-adjoint, we may assume that
  $$J_{22}=I_1\oplus -I_2\oplus0: \hbox{  } \hbox{  }\mathcal{M}_1\oplus\mathcal{M}_2\oplus\mathcal{M}_3,$$ where
  $\mathcal{M}_1=N(J_{22}-I),$  $\mathcal{M}_2=N(J_{22}+I)$ and  $\mathcal{M}_3=N(J_{22}).$
Obviously, $\mathcal{M}_1\oplus\mathcal{M}_2\oplus\mathcal{M}_3= R(P)^{\perp}.$
Suppose that $$J_{12}=\left(\begin{array}{ccc}
J_{12}^{(1)}&J_{12}^{(2)}&J_{12}^{(3)}\end{array}\right): \mathcal{M}_1\oplus\mathcal{M}_2\oplus\mathcal{M}_3\rightarrow R(P)$$  and
$$P_1=\left(\begin{array}{ccc}
P_{11}&P_{12} &P_{13}\end{array}\right): \mathcal{M}_1\oplus\mathcal{M}_2\oplus\mathcal{M}_3\rightarrow R(P).$$
So $J_{12}^{(1)}=0$ and $J_{12}^{(2)}=0$ follow from $J_{12}J_{22}=0.$
Moreover,  $J_{12}^*P_1+P_1^*J_{12}=J_{22}$ implies $$\left(\begin{array}{ccc}
0&0&P_{11}^*J_{12}^{(3)}\\0&0&P_{12}^*J_{12}^{(3)}\\ (J_{12}^{(3)})^*P_{11}&(J_{12}^{(3)})^*P_{12}& (J_{12}^{(3)})^*P_{13}+P_{13}^*J_{12}^{(3)}\end{array}\right)=I_1\oplus -I_2\oplus0,$$
 which induces $\mathcal{M}_1=\{0\}$ and $\mathcal{M}_2=\{0\}.$
 Then $J_{22}=0,$  so $J_{12}$ is a unitary operator with  $J_{12}^*P_1=-P_1^*J_{12}.$
 Setting $U:=\sqrt{-1}J_{12}^*,$  we conclude that $$J=\left(\begin{array}{cc}
0&\sqrt{-1}U^*
\\ -\sqrt{-1}U&0\end{array}\right):R(P)\oplus
R(P)^{\perp}$$ and $UP_1=P_1^*U^*.$
 Thus Lemma 3.2 implies $U=J_1V_1^*\oplus U_2.$\ \ \ $\Box$

In the following, the decomposition properties between symmetries $J\in\Delta_{P}$
 and symmetries $J'$ with $J'P\geq0$ are considered.
 It is clear that if $J\in\Delta_{P}$
  and the symmetry $J'$ satisfies $J'P\geq0,$  then
  $$JJ'PJ'J=JP^*J=I-P.$$
  Thus $J'':=JJ'\in\Gamma_{P}$ if and only if $JJ'=J'J.$
  However, we show that this situation $(JJ'=J'J)$ holds only if $P$ is an orthogonal projection.

 {\bf Lemma 3.6.} ([9, Corollary 14]) Let $P\in \mathcal{B(H)}^{Id}$ have the form (1.1). If $J$ is a symmetry, then
$JP\geqslant0$
if and only if $$J=\left(\begin{array}{cc}(I+P_1P_1^{\ast})^{-\frac{1}{2}}&(I+P_1P_1^{\ast})^{-\frac{1}{2}}P_1\\
P_1^{\ast}(I+P_1P_1^{\ast})^{-\frac{1}{2}}&J_2(I+P_1^{\ast}P_1)^{-\frac{1}{2}}\end{array}\right):R(P)\oplus R(P)^{\perp},$$
 where $J_2$ is a symmetry on the subspace
 $R(P)^{\perp}$  with $P_1=-P_1J_2.$

{\bf Lemma 3.7.} Let $P\in \mathcal{B(H)}^{Id}$  have the form (1.1). Then
\begin{equation}\begin{array}{rcl}\min\{J: JP\geq0,\ J=J^{\ast}=J^{-1}\}&=&2P_{(P+P^{\ast})^{+}}-I\\&=&(2P-I)|2P-I|^{-1}\\&=&\left(\begin{array}{cc}(I+P_1P_1^{\ast})^{-\frac{1}{2}}&(I+P_1P_1^{\ast})^{-\frac{1}{2}}P_1\\
P_1^{\ast}(I+P_1P_1^{\ast})^{-\frac{1}{2}}&-(I+P_1^{\ast}P_1)^{-\frac{1}{2}}\end{array}\right),\end{array} \end{equation}
where the ``min" is in the sense of Loewner partial order.

{\bf Proof.}  According to [9, Theorem 15] and [10, Remark],
we conclude that $$ min\{J: JP\geqslant0,\ J=J^{\ast}=J^{-1}\}=2P_{(P+P^{\ast})^{+}}-I=(P+P^{*}-I)|P+P^{*}-I|^{-1}$$ and
$$ 2P_{(P+P^{\ast})^{+}}-I=\left(\begin{array}{cc}(I+P_1P_1^{\ast})^{-\frac{1}{2}}&(I+P_1P_1^{\ast})^{-\frac{1}{2}}P_1\\
P_1^{\ast}(I+P_1P_1^{\ast})^{-\frac{1}{2}}&-(I+P_1^{\ast}P_1)^{-\frac{1}{2}}\end{array}\right):R(P)\oplus
R(P)^{\perp}.$$
  Then $$\min\{J: JP\geq0,\ J=J^{\ast}=J^{-1}\}=(2P-I)|2P-I|^{-1}=|2P-I|(2P-I)$$ follows from  Theorem 2.3.

{\bf Theorem 3.8.} Let $P\in \mathcal{B(H)}^{Id}$ and $J\in\Delta_{P}.$

(i) Then there exist uniquely symmetries $J_1$ and $J_2$ such that $J_1\in\Gamma_{P},$ $J_2P\geq 0$ and
$J=-\sqrt{-1}J_1J_2=\sqrt{-1}J_2J_1.$

(ii) If $P\notin \mathcal{P(H)},$ then there are not symmetries
$J_3$ and $J_4$ such that $J_3\in\Gamma_{P},$ $J_4P\geq 0$ and
$J=J_3J_4=J_4J_3.$

(iii) If $P\in \mathcal{P(H)},$ then there are symmetries
$J_3$ and $J_4$ such that $J_3\in\Gamma_{P},$ $J_4P\geq 0$ and
$J=J_3J_4=J_4J_3.$

{\bf Proof.} (i) Suppose that $P$ has the form (1.1). By Theorem 3.5, $JPJ=I-P^*$ yields $$J=\left(\begin{array}{cc}
0&\sqrt{-1}U^*
\\ -\sqrt{-1}U&0\end{array}\right):R(P)\oplus
R(P)^{\perp},$$ where $U\in \mathcal{B}(R(P), R(P)^{\perp})$ is a unitary operator with $UP_1=P_1^*U^*.$
Setting \begin{equation}J_2:=\left(\begin{array}{cc}(I+P_1P_1^{\ast})^{-\frac{1}{2}}&(I+P_1P_1^{\ast})^{-\frac{1}{2}}P_1\\
P_1^{\ast}(I+P_1P_1^{\ast})^{-\frac{1}{2}}&-(I+P_1^{\ast}P_1)^{-\frac{1}{2}}\end{array}\right):R(P)\oplus R(P)^{\perp},\end{equation} we conclude from Lemma 3.6 that $J_2$ is a symmetry with $J_2P\geq0.$
Define an operator $J_1$ as the form \begin{equation}J_1:=\left(\begin{array}{cc}
-U^*(I+P_{1}^{*}P_{1})^{-\frac{1}{2}}P_{1}^{*}&U^*(I+P_{1}^{*}P_{1})^{-\frac{1}{2}}
\\(I+P_{1}^{*}P_{1})^{-\frac{1}{2}}U&P_{1}^{*}U^{*}(I+P_{1}^{*}P_{1})^{-\frac{1}{2}}\end{array}\right):R(P)\oplus R(P)^{\perp}
.\end{equation} It is easy to see that $UP_1=P_1^*U^*$ implies
$$UP_1P_1^*=P_1^*U^*P_1^*=P_1^*P_1U,$$ so
$$U(I+P_1P_1^*)^{-\frac{1}{2}}=(I+P_1^*P_1)^{-\frac{1}{2}}U,$$ which induces $$(I+P_1P_1^*)^{-\frac{1}{2}}U^*=U^*(I+P_1^*P_1)^{-\frac{1}{2}}.$$
Then by a direct calculation, we get that  $J_1=J_1^*=J_1^{-1}$ and
$$\begin{array}{rl}J_1J_2=&\left(\begin{array}{cc}
0&-U^*P_1^*(I+P_1P_1^*)^{-1}P_1-U^*(I+P_1^*P_1)^{-1}
\\U(I+P_1P_1^*)^{-1}+U P_1P_1^*(I+P_1P_1^*)^{-1}&0\end{array}\right)
\\=&\left(\begin{array}{cc}
0&-U^*
\\ U&0\end{array}\right).\end{array}$$
Thus $J=-\sqrt{-1}J_1J_2=\sqrt{-1}J_2J_1.$

To show uniqueness, suppose that symmetries $J'_1$ and $J'_2$ satisfy that $J'_1PJ'_1=I-P,$ $J'_2P\geq 0$ and
$$J=-\sqrt{-1}J'_1J'_2=\sqrt{-1}J'_2J'_1.$$ Then
$J_1J_2=J'_1J'_2,$ which implies $J_1J'_1=J_2 J'_2.$
As $J'_1PJ'_1=I-P,$ we get from Theorem 2.7 that $J'_1$ has the following operator matrix form
\begin{equation}J'_1=\left(\begin{array}{cc}
-V(I+P_{1}^{*}P_{1})^{-\frac{1}{2}}P_{1}^{*}&V(I+P_{1}^{*}P_{1})^{-\frac{1}{2}}
\\(I+P_{1}^{*}P_{1})^{-\frac{1}{2}}V^{*}&P_{1}^{*}V(I+P_{1}^{*}P_{1})^{-\frac{1}{2}}\end{array}\right):R(P)\oplus R(P)^{\perp}
,\end{equation}
where $V$ is a unitary operator from $R(P)^{\bot}$ onto $R(P)$ with $V^{*}P_{1}=P_{1}^{*}V.$
Using Lemma 3.6 again, we know that $J'_2$ has the operator matrix form \begin{equation}J'_2=\left(\begin{array}{cc}(I+P_1P_1^{\ast})^{-\frac{1}{2}}&(I+P_1P_1^{\ast})^{-\frac{1}{2}}P_1\\
P_1^{\ast}(I+P_1P_1^{\ast})^{-\frac{1}{2}}&J'_{22}(I+P_1^{\ast}P_1)^{-\frac{1}{2}}\end{array}\right):R(P)\oplus R(P)^{\perp},\end{equation}
 where $J'_{22}$ is a symmetry on the subspace
 $R(P)^{\perp}$  with $P_1=-P_1J'_{22}.$
Combining equations (3.5), (3.6), $UP_1=P_1^*U^*$ and $V^{*}P_{1}=P_{1}^{*}V,$
 we conclude from a direct calculation that \begin{equation}J_1J'_1=U^*V^*\oplus U V.\end{equation}
 Also, equations (3.4) and (3.7) imply
 \begin{equation}J_2J'_2=\left(\begin{array}{cc}I&0\\
0&(I+P_1^{\ast}P_1)^{-1}P_1^{\ast}P_1-(I+P_1^{\ast}P_1)^{-1}J'_{22}\end{array}\right)
:R(P)\oplus R(P)^{\perp}.\end{equation}
 Then $J'_{22}=-I$ follows from equations (3.8), (3.9) and the fact $J_1J'_1=J_2 J'_2.$
Thus $J_2 =J'_2,$ so $J_1=J'_1.$

(ii) Conversely, we assume that there exist symmetries
$J_3$ and $J_4$ such that $J_3PJ_3=I-P,$ $J_4P\geq 0$ and
$J=J_3J_4=J_4J_3.$
Furthermore, Lemma 3.7 implies $J_4\geq J_2,$ which yields
$$\mathcal{P(H)}\ni P_4:=\frac{J_4+I}{2}\geq \frac{J_2+I}{2}:=P_2\in \mathcal{P(H)}.$$
Also, $J=J_3J_4=J_4J_3$ induces $J_3P_4=P_4J_3.$
Using (i), we get that $$-\sqrt{-1}J_1J_2=J=J_3J_4,$$ so $$
J_3J_4+\sqrt{-1}J_1J_2=0,$$ which implies $$J_3(2P_4-I)+\sqrt{-1}J_1(2P_2-I)=0.$$
 Therefore, $$P_4^\perp J_3P_4^\perp+\sqrt{-1}P_4^\perp J_1P_4^\perp=-P_4^\perp[J_3(2P_4-I)+\sqrt{-1}J_1(2P_2-I)]P_4^\perp=0,$$
where $P_4^\perp=I-P_4.$
Then $P_4^\perp J_3P_4^\perp=0$ and $ P_4^\perp J_1P_4^\perp=0$ follow from the fact that  $P_4^\perp J_3P_4^\perp$ and $P_4^\perp J_1P_4^\perp$ are self-adjoint.
  Thus $$J_3P_4^\perp=P_4^\perp J_3P_4^\perp=0,$$ which yields
   $P_4^\perp=0,$ that is $P_4=I.$ Hence $J_4=I$ and $J=J_3,$
   so $P=P^*.$ This is a contradiction with the fact $P\notin \mathcal{P(H)}.$

   (iii) is obvious.     \ \ \ $\Box$

 The following corollary is clear from Theorem 3.8 (ii).

 {\bf Corollary 3.9.}  Let $P\in \mathcal{B(H)}^{Id}$ and $J$ be a symmetry with $JP\geq0.$
If $P\notin \mathcal{P(H)},$  then there are not symmetries
$J_1\in\Gamma_{P}$ and $J_2\in\Delta_{P}$ such that
$J=J_1J_2=J_2J_1.$

{\bf Corollary 3.10.} Let $P\in \mathcal{B(H)}^{Id}$ and $J$ be a symmetry.
 Then $J\in\Delta_{P}$ if and only if $J'=\sqrt{-1}J(2P_{(P+P^{\ast})^{+}}-I)\in\Gamma_{P}.$

{\bf Proof.} Necessity is obvious from the proof of Theorem 3.8 (i) and Lemma 3.7.

Sufficiency. If $J'=\sqrt{-1}J(2P_{(P+P^{\ast})^{+}}-I)$ is a symmetry with
$J'PJ'=I-P,$  then $$J'(2P_{(P+P^{\ast})^{+}}-I)=-(2P_{(P+P^{\ast})^{+}}-I)J' \ \ \hbox{ and }\ \
 J=-\sqrt{-1}J'(2P_{(P+P^{\ast})^{+}}-I).$$
Also, Lemma 3.7 says $(2P_{(P+P^{\ast})^{+}}-I)P\geq0,$ so
$$(2P_{(P+P^{\ast})^{+}}-I)P(2P_{(P+P^{\ast})^{+}}-I)=P^*.$$
Thus $$JPJ=-J'(2P_{(P+P^{\ast})^{+}}-I)PJ'(2P_{(P+P^{\ast})^{+}}-I)=J'P^*J'=I-P^*.$$   \ \ \ $\Box$

{\bf Corollary 3.11.} Let $P\in \mathcal{B(H)}^{Id}$ and $J'$ be a symmetry with $J'P\geq0.$
If $dim N(P+P^{*})=dim N(2I-P-P^{*}),$ then the following statements are equivalent:

 (a) $J'=|2P-I|(2P-I);$

(b) $J'J=-JJ'$ for all symmetries $J$ with
$JPJ=I-P^*;$

(c) $J'J=-JJ'$ for a symmetry $J$ with
$JPJ=I-P^*.$

{\bf Proof.}  (a)$\Longrightarrow$(b) follows from Corollary 3.10. (b)$\Longrightarrow$(c) is obvious.

(c)$\Longrightarrow$(a). According to Theorem 3.5,
we know that $JPJ=I-P^*$ implies $$J=\left(\begin{array}{cc}
0&\sqrt{-1}U^*
\\ -\sqrt{-1}U&0\end{array}\right):R(P)\oplus
R(P)^{\perp},$$
where $U$ is a unitary operator from $R(P)$ onto $R(P)^{\perp}$ with  $UP_1=P_1^*U^*.$
Moreover, Lemma 3.6 yields $$J'=\left(\begin{array}{cc}(I+P_1P_1^{\ast})^{-\frac{1}{2}}&(I+P_1P_1^{\ast})^{-\frac{1}{2}}P_1\\
P_1^{\ast}(I+P_1P_1^{\ast})^{-\frac{1}{2}}&J_2(I+P_1^{\ast}P_1)^{-\frac{1}{2}}\end{array}\right):R(P)\oplus R(P)^{\perp},$$
 where $J_2$ is a symmetry on the subspace
 $R(P)^{\perp}$  with $P_1=-P_1J_2.$
 Thus \begin{equation}\sqrt{-1}(I+P_1P_1^{\ast})^{-\frac{1}{2}}U^*=-\sqrt{-1}U^*J_2(I+P_1^{\ast}P_1)^{-\frac{1}{2}} \end{equation} follows from equation $J'J=-JJ'.$
 As $UP_1=P_1^*U^*$ implies $$ (I+P_1P_1^{\ast})^{-\frac{1}{2}}U^*=U^*(I+P_1^{\ast}P_1)^{-\frac{1}{2}},$$
we conclude from equation (3.10) that $J_2=-I.$ Then Lemma 3.7 yields
 $J'=|2P-I|(2P-I).$
 \ \ \ $\Box$

 Combining Corollary 3.10 and 3.11, we get the following corollary.

 {\bf Corollary 3.12.}  Let $P\in \mathcal{B(H)}^{Id}$ with  $dim N(P+P^{*})=dim N(2I-P-P^{*}).$
 Then $$\begin{array}{rl}&\{J: JP\geq0,\ J=J^{\ast}=J^{-1}\}
 \cap \{J: JJ'=-J'J \hbox{ for a  } J'\in\Delta_{P} \}\\=&\{|2P-I|(2P-I)\}\\=&\{J: JP\geq0,\ J=J^{\ast}=J^{-1}\}
 \cap \{J: JJ'=-J'J \hbox{ for a } J'\in\Gamma_{P}\}.\end{array}$$

\end{document}